\documentclass{siamart190516}
\usepackage{epsfig,amsmath,amsfonts,amssymb,color}
\usepackage{algorithmicx,algpseudocode}

\newcommand{\dist}{\operatorname{dist}}

\newcommand{\xb}{\mathbf{x}}
\newcommand{\yb}{\mathbf{y}}
\newcommand{\rb}{\mathbf{r}}
\newcommand{\hb}{\mathbf{h}}
\newcommand{\R}{\mathbb{R}}
\newcommand{\C}{\mathbb{C}}

\newcommand{\V}{\mathcal{V}}
\newcommand{\K}{\mathcal{K}}
\newcommand{\N}{\mathcal{N}}

\newcommand{\Nbrs}{\mathcal{N}}
\newcommand{\Inter}{\mathcal{I}}
\newcommand{\Kids}{\mathcal{K}}
\newcommand{\Triang}{\mathcal{T}}
\newcommand{\TriangS}{\mathcal{S}}
\newcommand{\Ord}{\mathcal{O}}
\newcommand{\half}{\frac{1}{2}}
\newcommand{\thrf}{\frac{3}{2}}
\newcommand{\abs}[1]{\left| #1 \right|}
\newcommand{\norm}[2]{{\|#1\|}_{#2}}
\newcommand{\inner}[3]{\left\langle #1 , #2 \right\rangle_{#3}}

\begin{document}

\title{A fast method for evaluating Volume potentials in the Galerkin
  boundary element method
\thanks{This work was in part funded by the National Science
  foundation under the grant DMS-1720431.}}
\author{Sasan Mohyaddin \and Johannes Tausch\thanks{Department of Mathematics, 
Southern Methodist University, Dallas, TX 75275,
\email{smohyaddin@smu.edu,tausch@smu.edu}}}

\maketitle

\begin{abstract}
  Three algorithm are proposed to evaluate volume potentials that
  arise in boundary element methods for elliptic PDEs.  The approach
  is to apply a modified fast multipole method for a boundary
  concentrated volume mesh. If $h$ is the meshwidth of the boundary,
  then the volume is discretized using nearly $\Ord(h^{-2})$ degrees
  of freedom, and the algorithm computes potentials in nearly
  $\Ord(h^{-2})$ complexity. Here nearly means that logarithmic terms
  of $h$ may appear. Thus the complexity of volume potentials
  calculations is of the same asymptotic order as boundary potentials.
  For sources and potentials with sufficient regularity the parameters of the
  algorithm can be designed such that the error of the approximated
  potential converges at any specified rate $\Ord(h^p)$.  The accuracy
  and effectiveness of the proposed algorithms are demonstrated for
  potentials of the Poisson equation in three dimensions.
\end{abstract}

\begin{keywords}
Fast multipole method, boundary integral equation, volume potential,
boundary concentrated mesh, Poisson equation.
\end{keywords}

\begin{AMS}
65N38, 65N12, 65N30
\end{AMS}

\section{Introduction}\label{sec:intro}
Boundary integral methods for homogeneous elliptic partial
differential equations are based on representing the solution in form
of layer potentials. This results in an integral equation on the
boundary of the domain. For a three dimensional domain this implies a
reduction to a problem on the two dimensional boundary. After
discretization with a boundary mesh of size $h$, one obtains a dense
matrix of size $O(h^{-2})$. Typically, the linear system is solved
iteratively, where the dominant numerical cost is the evaluation of
matrix-vector products. There are several well established methods to
accelerate this operation. This includes the fast multipole method
\cite{greengard-rokhlin87,na-ko-le-wh94,tausch04}, wavelets
\cite{dahmen-harbrecht-schneider07} and $\mathcal{H}$-matrix algebra
\cite{borm-grase-hackb02} which can be combined with adaptive cross
approximation~\cite{bebendorf08}. With these methods it is possible to
approximately compute the matrix-vector product with nearly or even
exactly $O(h^{-2})$ complexity, while maintaining the convergence rate
of the discretization error.

If the underlying PDE is inhomogeneous, the solution must be
represented with an additional volume potential of the right hand side
of the equation. Likewise, the reconstruction of the solution in the
domain requires the evaluation of layer potentials in the volume.  The
efficient evaluation volume potentials has been the subject of many
investigations.  A popular method is the dual reciprocity method. Here
the basic idea is to approximate the right hand side by radial basis
functions, and to use integration by parts to convert the volume
integral to a boundary integral, see, e.g.,
\cite{partr-breb-wrob}. The approach in \cite{etheridge-greengard01}
is based on related ideas.  Another frequently used approach is to
embed the domain into a rectangular box and apply either the fast
multipole method \cite{mcKenny-gg95,askham-cerfon17} or a fast Poisson
solver in the box. To avoid difficulties extending the right hand side
beyond the domain one can discretize the volume, for instance, with a
tetrahedral mesh and use a fast method for the evaluation of the
domain integral, see \cite{of-stein-urthaler10}.

We also mention some methods for two dimensional domains that either
rely on a Fourier-Galerkin discretization of the boundary curve
\cite{guan-etal21} or are specific to circular
domains~\cite{atkinson85,steinbach-tchoualag14}. A comparison of
different domain evaluation techniques is given in
\cite{ingber-mammoli-brown01}.

However, the order of the complexity will be increased when volume
potentials appear in the integral equation.  Likewise, the evaluation
of the solution will increase the complexity if a uniform volume
mesh is used for the domain evaluations of layer or volume potentials.
We will therefore consider discretizations using a boundary
concentrated (BC) mesh, where the meshwidth of the volume
discretization grows proportionally with the distance from the
boundary.  This kind of mesh has already been employed in the context
of finite element methods \cite{khoromskij-melenk03}.  The number
elements in such a mesh is order $O(h^{-2})$ where $h$ is the
meshwidth of the boundary mesh. Hence the number of elements of in the
boundary and volume meshes have the same asymptotic order. In this
article, we will derive fast algorithms for the following computational
tasks.
\begin{itemize}
\item{\textbf{Volume to volume (VtV).}}
  Given a function represented by the BC mesh, compute its volume potential
  on the BC mesh.
\item{\textbf{Volume to boundary (VtB).}} Given a function represented by the BC
  mesh, compute the volume potential on the surface mesh.
\item{\textbf{Boundary to volume (BtV).}} Given a density represented by the surface mesh,
  compute its layer potential on the BC mesh.
\end{itemize}
In particular, we will devise a fast multipole algorithm for BC meshes and
show that its parameters can be chosen such that it has nearly optimal  
$O(h^{-2})$ complexity.

This paper is organized as follows. In the remainder of this section
we provide more detailed background material on layer and volume
integrals and their discretization. A hierarchical subdivision of the
volume by a boundary concentrated meshes is then described in
section~\ref{sec:triangulation}. Section~\ref{sec:bcfmm} describes a
fast multipole type algorithm to efficiently perform
the VtV, VtB and BtV calculations. Section~\ref{sec:analysis:err}
provides an analysis of the complexity and
accuracy of the methods. We conclude in section~\ref{sec:numresults}
with numerical results that illustrate the theoretical results.

\subsection{Boundary and Volume Potentials}
Consider an elliptic operator 
$\mathcal{L}$ with constant coefficients
in a bounded
domain $\Omega \subset \R^3$ with boundary surface $\Gamma = \partial \Omega$.
The Green's representation formula expresses the solution of
$\mathcal{L}u =f$
in $\Omega$ 
in terms of the Dirichlet and Neumann data on $\Gamma$
\begin{equation}\label{green:ie:vol}
  u = \tilde \V [\gamma_1 u] - \tilde \K [\gamma_0 u] + \tilde \N f,
  \qquad \text{in}\;\Omega.
\end{equation}
Here, $\gamma_0 u$ and $\gamma_1 u$ are the boundary trace and normal
boundary trace of a function $u$ defined in the domain, and
$\tilde \V$ and $\tilde \K$ are the single-layer and double-layer
potentials, defined by
\begin{align*}
  \tilde \V q(\xb) &= \int_{\Gamma} G(\xb,\yb) q(\yb)\, ds_\yb\,,\\
  \tilde \K u(\xb) &= \int_{\Gamma} \frac{\partial G}{\partial n_\yb}(\xb,\yb) u(\yb)\, ds_\yb.
\end{align*}
Moreover, $\tilde \N$ denotes the volume (or sometimes Newton) potential
\begin{equation*}
  \tilde \N f(\xb) = \int_{\Omega} G(\xb,\yb) f(\yb)\, d\yb.
\end{equation*}
The kernel $G(\cdot,\cdot)$ is the free space Green's function of
the PDE, which in the case of the Poisson equation is
\begin{equation*}
  G(\xb,\yb) = \frac{1}{4\pi} \frac{1}{\abs{\xb-\yb}}.
\end{equation*}

Taking the boundary trace in \eqref{green:ie:vol} results in the Green's
integral equation
\begin{equation}\label{green:ie:bdry}
  \half \gamma_0 u = \V [\gamma_1 u] - \K [\gamma_0 u] + \N f,
  \qquad \text{on}\;\Gamma.
\end{equation}
In \eqref{green:ie:vol} the operators with a tilde
indicate that a potential is evaluated in the domain, and in \eqref{green:ie:bdry}
the potentials without the tilde indicates the evaluation on the
boundary. It is well known that the potentials are continuous in the following
spaces
\begin{align*}
  \tilde \V : H^{-\half}(\Gamma) \to H^1(\Omega), & \quad
         \V:H^{-\half}(\Gamma) \to H^{\half}(\Gamma),\\
  \tilde \K : H^{\half}(\Gamma) \to H^1(\Omega) ,& \quad
         \K:H^{\half}(\Gamma) \to H^{\half}(\Gamma),\\
  \tilde \N : H^{-1}(\Omega) \to H^1(\Omega), & \quad
         \N:H^{-1}(\Omega) \to H^{\half}(\Gamma),
\end{align*}
see, e.g., \cite{hsiao-wendland08,mclean00}.
In the direct boundary element
method, the integral equation \eqref{green:ie:bdry} is solved for the
missing boundary data, and then the solution in the interior is
evaluated using the Green's representation formula.

\subsection{Discretization}
We briefly describe the Galerkin discretization of surface and volume
integral operators.  To fix ideas, assume that $\Omega$ is a
polyhedral domain that has been subdivided into a tetrahedral mesh
$\Triang$. We assume that this volume mesh is shape regular, but not
necessarily conforming or quasi-uniform. However, we assume that the
restriction  to the boundary is a conforming, shape regular, and
quasi-uniform triangular mesh of meshwidth $h$.  The boundary mesh is
denoted by $\TriangS$. 

Suppose we want to solve the Dirichlet problem using the direct
integral formulation. In this case the Green's integral formulation
\eqref{green:ie:bdry} is solved for the Neumann data, and then the
representation formula \eqref{green:ie:vol} provides the solution in
the domain.  

The Galerkin discretization of the integral equation is based on the
finite element space $S^\Gamma_h$. A typical choice for $S^\Gamma_h$
are low-order piecewise polynomial functions on $\TriangS$. The
Galerkin discretization of \eqref{green:ie:bdry} reads: find
$q_h \in S_h^\Gamma$ such that
\begin{equation}\label{galerkin:ie}
  \inner{\varphi}{\V q_h}{L_2(\Gamma)}
  = \inner{\varphi}{\left(1/2 + \K\right)g}{L_2(\Gamma)}
   - \inner{\varphi}{\N f}{L_2(\Gamma)}
\end{equation}
for all basis functions $\varphi$ of $S^\Gamma_h$. 
Here $\inner{\cdot}{\cdot}{L_2(\Gamma)}$ is the $L_2$-inner product on $\Gamma$, 
$q_h$ is the approximation of the Neumann data and $g$ is the given
Dirichlet data. Since $q_h$ and is a linear
combination of all $\varphi$'s this leads to a linear system
\begin{equation*}
  V \mathbf{q} = \mathbf{b},
\end{equation*}
where $\mathbf{q}$ is the vector of coefficients in the expansion of
$q_h$ in the basis, and the coefficients of $\mathbf{b}$ are given by
the right hand side in \eqref{galerkin:ie}. 

For the approximation of the solution in the volume we also use a
variational approach. To that end, the potential $u$ is approximated
in the space of piecewise polynomial functions $S^\Omega_h$ on
$\Triang$. Since in BC meshes the size of the elements varies, the
polynomial order $p_\omega$ is tied to the size of the elements. The
precise relationship will be discussed in
section~\ref{sec:analysis:err} below.

For a tetrahedron orthogonal polynomials
$\phi^\alpha_\omega, \abs{\alpha}\leq p_\omega$ can be constructed
explicitly in terms of Jacobi polynomials, see \cite{koornwinder75}. Outside of
$\omega$ these functions are extended by zero.
Thus the finite element space is
\begin{equation}\label{def:femspace:volume}
S^\Omega_h = \mbox{span}\left\{ \phi^\alpha_\omega \;:\; \omega \in
  \mathcal{T},\; \abs{\alpha} \leq p_\omega\right\}.
\end{equation}
The space $S^\Omega_h$ is a subset of $L_2(\Omega)$, but it is not contained in
$H^1(\Omega)$. However, this is sufficient regularity for the
$L_2(\Omega)$-orthogonal projection which is given as follows,
\begin{align}
  P_\Triang u &= \sum_{\omega \in \mathcal{T}} \sum_{\abs{\alpha} \leq p_\omega}
  u_{\alpha,\omega}{\phi^\alpha_\omega}, \label{galerkin:vol}\\
  \intertext{where}
  u_{\alpha,\omega}&= \inner{\phi^\alpha_\omega}
   {\tilde \V [\gamma_1 u] - \tilde \K [\gamma_0 u] + \tilde \N f}{L_2(\Omega)}.
\end{align}
In analogy to $S^\Omega_h$, the space $S^\Gamma_h$ is spanned by
piecewise polynomials on triangles. For the latter all triangles have
diameter proportional to $h$ and the polynomial
degree $p$ is fixed and typically low. We write
\begin{equation}\label{def:femspace:surface}
S^\Gamma_h = \mbox{span}\left\{ \varphi^\alpha_\gamma \;:\; \gamma \in
  \mathcal{S},\; \abs{\alpha} \leq p \right\}.
\end{equation}
Here $\alpha$ is a multi index with two components whereas in
\eqref{def:femspace:volume} $\alpha$ has three components.

\section{Boundary Concentrated Mesh}\label{sec:triangulation}
This section provides more details for the case that $\Triang$ is a BC
mesh.  We consider a polyhedron $\Omega \subset \R^3$ that is
subdivided into a small number of tetrahedra. These tetrahedra form
the coarsest level, or level $\ell=0$ in a hierarchical subdivision of
space. The $\ell+1$-st level tetrahedra are obtained by subdividing
some or all tetrahedra in the $\ell$-th level into eight
tetrahedra. The refinements of a tetrahedron $\omega$ are referred to
as the children, or $\Kids(\omega)$.

Note that some care must be applied to ensure that the
refinements remain shape regular, because in general,
it is not possible to obtain congruent subdivisions of three
dimensional tetrahedra. However, with the refinement scheme introduced
in \cite{bey95} one can limit the number of congruency classes to
three which implies shape regularity.

A quasi uniform refinement is achieved if all tetrahedra in a given
level are subdivided.  In this case a tree structure results, where
the children of the root are the tetrahedra of the initial subdivision
of $\Omega$. All other nodes either have eight children or are leaves
in the finest level.  The domain $\Omega$ is the union of all
tetrahedra in any given level.

\begin{figure}
\begin{center}
\includegraphics[height=4cm]{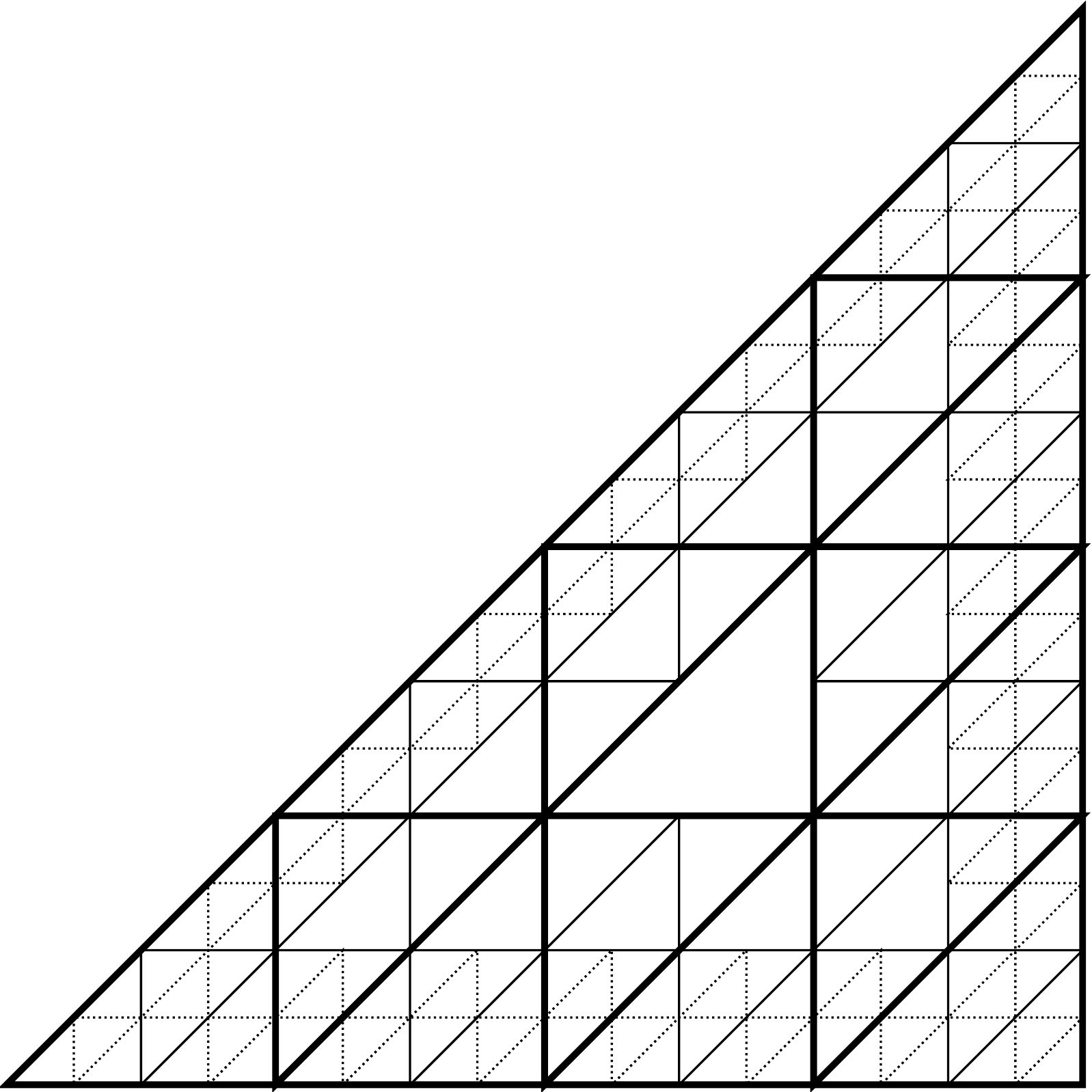}
\end{center}
\caption{
  A two dimensional illustration of a boundary concentrated refinement.} 
\label{fig:bcRefinement}
\end{figure}

In contrast to uniform refinements, a boundary concentrated refinement
is obtained by refining only tetrahedra that are close to the
boundary. A two dimensional situation is illustrated in
figure~\ref{fig:bcRefinement}. The resulting tree structure may have
leaves in any level, which are the tetrahedra that have not been
refined. The domain $\Omega$ is the union of all leaves in all levels.

To describe the BC refinement scheme in
detail, denote the center of $\omega$
by $\xb_\omega$, the diameter by
\begin{equation}\label{def:diam}
\rho_\omega = \max_{{\bf v}: \mbox{\scriptsize vertex of }\omega} \abs{{\bf v}-\xb_\omega}.
\end{equation}
The separation ratio of two tetrahedra in the same level is
defined as
\begin{equation}\label{def:sepRatio}
\eta(\omega,\omega') = \frac{\rho_{\omega'} + \rho_{\omega'}}{\abs{\xb_\omega - \xb_{\omega'}}}\,,
\quad \omega \not= \omega',
\end{equation}
and  $\eta(\omega,\omega) = \infty$.
The neighbors of $\omega$ are the tetrahedra $\omega'$ in the same level 
for which the separation ratio is greater than a predetermined
constant  $\eta_0$. That is,
\begin{equation}\label{def:neighbors}
  \Nbrs(\omega) = \{ \omega' \in C_\ell \,:\, \eta(\omega,\omega') > \eta_0  \}.
\end{equation}
Here, $C_\ell$, denotes the set of all tetrahedra in level $\ell$.
Further,
\begin{equation*}
B_\ell =
\{ \omega \in C_\ell \,:\, \omega\; \text{has a face in}\;\Gamma \}
\end{equation*}
denotes the set of all boundary tetrahedra in level $\ell$. Here it is worth
emphasizing that if $\omega$ only has one vertex or one edge in $\Gamma$ it
is not included in $B_\ell$. Moreover,
\begin{equation*}
M_\ell =
\{ \omega \in C_\ell \,:\, \Nbrs(\omega) \cap B_\ell \not= \emptyset\}
\end{equation*}
denotes the tetrahedra that have a boundary tetrahedron among their
neighbors. These are the tetrahedra are marked for refinement. Thus
\begin{equation*}
L_\ell = C_\ell \setminus M_\ell\,,
\end{equation*}
are the leaves in level $\ell$. The next
level list of tetrahedra is
\begin{equation*}
  C_{\ell+1} = \bigcup\limits_{\omega \in M_\ell}\Kids(\omega).
\end{equation*}
The refinement process is repeated until a finest level $L$ is reached. There,
all tetrahedra are leaves, but we distinguish between tetrahedra
near the surface and tetrahedra away from the surface. Hence we set
$L_L = C_L \setminus M_L$ and denote by $L_\ell^*$ the set of leaves
in any level, i.e.,
\begin{equation*}
  L_\ell^* = L_\ell,\; \ell\in \{0,\dots L-1\} \quad\mbox{and}\quad L_L^* = C_L.
\end{equation*}
We obtain the following subdivision of $\Omega$
\begin{equation}\label{spaceDecomp:Om}
  \Triang =
  M_L \cup \bigcup\limits_{\ell=0}^{L} L_\ell
  =
  \bigcup\limits_{\ell=0}^L L_\ell^*\,.
\end{equation}
An important concept in the fast multipole method is the interaction
list $\Inter(\omega)$, which consists of tetrahedra whose parents are
neighbors of the parent of $\omega$, but who are not neighbors with
$\omega$ itself. In level zero, we set
$\Inter(\omega) = C_0 \setminus \Nbrs(\omega)$, which can possibly be
an empty set.

Since $\Triang$ is a geometric mesh, it is known that its cardinality
is order $4^L \sim h^{-2}$, where the constant
depends on $\eta_0$. It remains to ensure
is that the number of neighbors and interacting tetrahedra is
bounded. 

\begin{theorem}\label{theo:cplx}
The cardinalities $\#\Nbrs(\omega)$ and $\#\Inter(\omega)$ are uniformly
bounded. Furthermore, there is a constant $c$ such that
\begin{equation*}
  \#B_\ell \leq c 4^{\ell}, \quad  
  \#M_\ell \leq c 4^{\ell}, \quad  
  \#L_\ell \leq c 4^{\ell}, \quad  \mbox{and } \quad
  \#C_\ell \leq c 4^{\ell}.
\end{equation*}
\end{theorem}

\begin{proof}
Two tetrahedra $\omega$, $\omega'$ in $C_\ell$ are neighbors if
$\abs{\xb_\omega - \xb_{\omega'}} < (\rho_\omega + \rho_{\omega'})/\eta_0$.
If we let $R_\ell = \max\{ \rho_\omega : \omega \in C_\ell\}$ then it
follows that all neighbors of $\omega$ are contained in the sphere
$B_\omega$ with center $\xb_\omega$ and
radius $(1+2/\eta_0)R_\ell$.

Further, let $\tilde \rho_\omega$ be the radius of the largest sphere
that is contained in $\omega$, and let $\tilde R_\ell = \min\{ \tilde
  \rho_\omega : \omega \in C_\ell\}$. Shape regularity implies that 
$R_\ell/\tilde R_\ell \leq c$. Since the enclosed spheres of the
neighbors are contained in $B_\omega$ it follows for their volumes that
\begin{equation*}
\#\Nbrs(\omega) \frac{4\pi}{3} \tilde R_\ell^3
\leq
\frac{4\pi}{3} \sum_{\omega'\in\Nbrs(\omega)} \tilde \rho_{\omega'}^3
\leq 
|B_\omega|
=
\frac{4\pi}{3}  \left(1+\frac{2}{\eta_0}\right)^3 R_\ell^3,
\end{equation*}
which implies that
\begin{equation*}
\#\Nbrs(\omega) 
\leq
\frac{R_\ell^3}{\tilde R_\ell^3 }  \left(1+\frac{2}{\eta_0}\right)^3 
\end{equation*}
so the number of neighbors is indeed bounded. The boundedness of $
\#\Inter(\omega)$ follows immediately from the definition of
interaction lists. 

Since all boundary faces are refined into four faces in each step and
since each boundary
face belongs to only one boundary tetrahedron it follows that
$\#B_\ell \leq  c 4^{\ell}$. Further, since every $\omega\in M_\ell$
is a neighbor of a boundary tetrahedron it follows that $\#M_\ell \leq
c 4^{\ell}$. The
remaining estimates follow from $\#L_\ell \leq \#C_\ell = 8  \#M_{\ell-1}$. 
\end{proof}

\section{Boundary Concentrated FMM}\label{sec:bcfmm}
The fast multipole method is based on a hierarchical splitting
of the source and target domains. In the standard method, this hierarchy
can be viewed as a tree with all leaves in the finest level. On the
other hand, the
BC refinement leads to a tree with leaves in any level. This implies some
modifications for the calculations for the nearfield which we describe
in this section. 

A key idea is to break neighbor interactions in
a given level into neighbors and farfield interactions in the next
finer levels. Suppose for now that level $\ell$ has no leaf nodes,
then this can be written as
\begin{equation}\label{splitlevelnoleaves}
L_\ell = \emptyset \; \Rightarrow
  \bigcup\limits_{\omega \in C_\ell} \omega\times \Nbrs(\omega)
  = \bigcup\limits_{\omega \in C_{\ell+1}} \omega\times \Nbrs(\omega)
  \;\; \cup \;
   \bigcup\limits_{\omega \in C_{\ell+1}} \omega\times \Inter(\omega).
\end{equation}
Here the left and right sets in a Cartesian product indicate targets
(i.e., $\xb$-variable) and sources (i.e., $\yb$-variable) of the
volume potential operator.
 
If there are leaves, this splitting gets more complicated, since some nodes
have refinements in the next level, whereas others do not. The
neighbors of any $\omega \in C_\ell$ may contain leaves and marked
nodes. To distinguish them we set
\begin{equation*}
  \Nbrs_L(\omega) = \Nbrs(\omega) \cap L_\ell,
  \quad\text{and}\quad
  \Nbrs_M(\omega) = \Nbrs(\omega) \cap M_\ell\,.
\end{equation*}
Since leaf neighbors have no refinements, it turns out that the neighbor lists have to
be extended to contain nodes from different levels. To
that end, denote by $\ell(\omega)$ the level of $\omega$
and by $\pi_\ell(\omega)$ the parent of $\omega$ in level $\ell$. Then
the extended neighbor list of $\omega$ is defined as
\begin{equation*}
  \Nbrs^*(\omega) = \Nbrs_M(\omega) \;\cup\;
  \bigcup_{\ell=0}^{\ell(\omega)} \Nbrs_L(\pi_\ell(\omega)).
\end{equation*}
Here $\pi_\ell(\omega) = \omega$ if $\ell = \ell(\omega)$. 
The following lemma generalizes \eqref{splitlevelnoleaves} for the
case that leaves are present in a given refinement level.

\begin{lemma}\label{lem:splitnear}
  It holds that
 \begin{equation*}
   \bigcup_{\omega \in M_\ell} \omega \times \Nbrs^*(\omega) =
   \bigcup_{\omega \in M_{\ell+1}} \omega \times \Nbrs^*(\omega)
   \;\cup \bigcup_{\omega \in L_{\ell+1}} \omega \times \Nbrs^*(\omega)
   \;\cup \bigcup_{\omega \in C_{\ell+1}} \omega \times \Inter(\omega).
\end{equation*}
\end{lemma}
From theorem~\ref{theo:cplx} it follows that 
$\#\Nbrs^*(\omega) \leq c L$, thus the extended nearfield computations
will contribute only a logarithmic term in the complexity of
algorithm~\ref{algo:vtv}.

\begin{proof}[Proof of Lemma \ref{lem:splitnear}]
The definition of extended neighbors implies that
\begin{equation*}
  \bigcup_{\omega \in M_\ell} \omega \times \Nbrs^*(\omega) =
  \bigcup_{\omega \in M_\ell} \omega \times \Nbrs_M(\omega) \;\cup
  \bigcup_{\omega \in M_{\ell} \atop \ell'=0..\ell } \omega \times \Nbrs_L(\pi_{\ell'}(\omega)).
\end{equation*}
In the first term both sources and targets have refinements, hence
a splitting into neighbors and interaction lists of the next level as
in \eqref{splitlevelnoleaves} can be performed. In the second term only
$\omega$ can be refined. This leads to
\begin{align*}
  \bigcup_{\omega \in M_\ell} \omega \times \Nbrs^*(\omega) &=
  \bigcup_{\omega \in C_{\ell+1}} \omega \times \Nbrs(\omega) 
  \;\cup \bigcup_{\omega \in C_{\ell+1} \atop \ell'=0..\ell } \omega \times \Nbrs_L(\pi_{\ell'}(\omega))
  \;\cup \bigcup_{\omega \in C_{\ell+1}} \omega \times \Inter(\omega)\\
  &= \bigcup_{\omega \in C_{\ell+1}} \omega \times \Nbrs_M(\omega)
  \;\cup \bigcup_{\omega \in C_{\ell+1} \atop \ell'=0..\ell+1 } \omega \times \Nbrs_L(\pi_{\ell'}(\omega))
  \;\cup \bigcup_{\omega \in C_{\ell+1}} \omega \times \Inter(\omega)\\
  &= \bigcup_{\omega \in C_{\ell+1}} \omega \times \Nbrs^*(\omega)
  \;\cup \bigcup_{\omega \in C_{\ell+1}} \omega \times \Inter(\omega)
\end{align*}
In the second step the leaf neighbors in the first
term are incorporated into the second term. The third step follows from the definition
of extended neighbors. Splitting the first term of the last equation
into $C_{\ell+1}= M_{\ell+1} \cup L_{\ell+1}$ gives the assertion.
\end{proof}

\subsection{Decomposition for the VtV calculation}
We now turn to the decomposition of the source and target domains into
nearfields and farfields, which rely on a repeated application of
Lemma~\ref{lem:splitnear}. In the coarsest level, 
$\Nbrs(\omega) = \Nbrs^*(\omega)$ and $\Omega = \Nbrs^*(\omega)
\cup \Inter(\omega)$ holds for $\omega \in C_0$. Further, $\Omega = C_0
= M_0 \cup L_0$.  Thus 
\begin{equation*}
  \Omega \times \Omega = 
  \bigcup\limits_{\omega \in M_0} \omega \times \Nbrs^*(\omega)\;\; \cup \;
  \bigcup\limits_{\omega \in L_0} \omega \times \Nbrs^*(\omega)\;\; \cup \;
  \bigcup\limits_{\omega \in C_0} \omega \times \Inter(\omega)
\end{equation*}
Applying lemma~\ref{lem:splitnear} to the first term gives
\begin{equation*}
  \Omega \times \Omega = 
  \bigcup\limits_{\omega \in M_1} \omega \times \Nbrs^*(\omega)\;\; \cup \;
  \bigcup\limits_{\omega \in L_0\cup L_1} \omega \times \Nbrs^*(\omega)\;\; \cup \;
  \bigcup\limits_{\omega \in C_0\cup C_1} \omega \times \Inter(\omega)
\end{equation*}
holds. Hence it follows by recursion through levels that
\begin{equation}\label{spaceDecomp:OmOm}
  \Omega \times \Omega =
 \bigcup\limits_{\ell=0}^L  \bigcup\limits_{\omega \in L^*_\ell} 
    \omega \times \Nbrs^*(\omega) \;\;\cup \;\;
 \bigcup\limits_{\ell=0}^L  \bigcup\limits_{\omega \in C_\ell} 
    \omega \times \Inter(\omega)
\end{equation}
In this decomposition the nearfield involves nodes in all
levels and sources and targets can be in different levels.
However, the cardinality of the sets $C_\ell$ is much smaller with
a boundary concentrated refinement than with a uniform refinement.
Therefore the algorithm based on \eqref{spaceDecomp:OmOm} will be more
efficient.

\subsection{Decomposition for the VtB, BtV calculations}
For these calculations either the source or the target domain are replaced
by the boundary. We denote the set of boundary faces of a tetrahedron
by $\gamma(\omega)$, further
\begin{equation*}
  \Nbrs_{\Gamma}(\omega) \;=
  \bigcup\limits_{\omega' \in \Nbrs(\omega)}  \gamma(\omega')
\quad\text{and}\quad
  \Inter_{\Gamma}(\omega) \;=
  \bigcup\limits_{\omega' \in \Inter(\omega)}  \gamma(\omega').
\end{equation*}
By the definition of the set $L_\ell$ it follows that $\gamma(\omega)
= \emptyset$ and  $\Nbrs_{\Gamma}(\omega) = \emptyset$ when $\omega
\in L_\ell$. Thus $\Nbrs_{\Gamma}(\omega)$ is the restriction of
$\Nbrs^*(\omega)$ to $\Gamma$.
Moreover, $\gamma(\omega)$ is non-empty if and only if  
$\gamma(\omega)$ is in $B_\ell$. With this in mind the appropriate
domain decompositions
can be obtained from \eqref{spaceDecomp:OmOm} by restricting either
the source or target domain to the boundary. It follows that
  \begin{align}
  \Gamma \times \Omega &=
  \bigcup\limits_{\omega \in B_L} \gamma(\omega) \times \Nbrs(\omega)\;\; \cup \;\;
 \bigcup\limits_{\ell=0}^L  \bigcup\limits_{\omega \in B_\ell} 
                         \gamma(\omega) \times \Inter(\omega),\label{spaceDecomp:GaOm}\\
  \Omega \times \Gamma &=
  \bigcup\limits_{\omega \in L^*_L} \omega \times \Nbrs_\Gamma(\omega)\;\; \cup \;\;
 \bigcup\limits_{\ell=0}^L  \bigcup\limits_{\omega \in C_\ell} 
    \omega \times \Inter_\Gamma(\omega)\,.\label{spaceDecomp:OmGa}
\end{align}
These decompositions are considerably simpler than
\eqref{spaceDecomp:OmOm}, because the nearfields only involve finest
level tetrahedra and sources and targets are always in the same level.

\subsection{Translation operators}
The evaluation of the farfield can be accomplished with the
moment-to-local translation of the fast multipole algorithm. For
completeness, we briefly recall its derivation for the case that the
kernel is approximated by a truncated Taylor series expansion. More
details can be found, e.g., in \cite{tausch04}.

The volume potential of source $\omega'$ is denoted by
$\tilde \N_{\omega'} f(\xb)$. When $\omega' \in \Inter(\omega)$ and $\xb \in \omega$
then this potential can be approximated by the $q_\ell$-th order Taylor expansion
centered at $(\xb_\omega, \xb_{\omega'})$. Thus
\begin{align*}
\tilde \N_{\omega'} f(\xb) &= \int_{\omega'} 
  G(\xb, \yb) f(\yb) \, d\yb \\
&\approx
\sum_{\abs{\alpha} \leq q_\ell} \sum_{\abs{\beta} \leq q_\ell - \abs{\alpha}} 
\frac{D^{\alpha+\beta} G(\xb_\omega, \xb_{\omega'})}{\alpha!\beta!}
(\xb - \xb_\omega)^\alpha 
\int_{\omega'}(\xb_{\omega'} - \yb)^\beta f(\yb)\, d\yb\\
&=
\sum_{\abs{\alpha} \leq q_\ell}
\lambda_{\omega}^\alpha (\xb - \xb_\omega)^\alpha.
\end{align*}
In the formula above, the expansion coefficients
$\lambda_{\omega}^\alpha$ are given by 
\begin{equation*}
\lambda_{\omega}^\alpha = \sum_{\abs{\beta} \leq q_\ell - \abs{\alpha}}
\frac{D^{\alpha+\beta} G(\xb_\omega, \xb_{\omega'})}{\alpha!\beta!}
(-1)^{\abs{\beta}} m^\beta_{\omega'}\!(f),\quad \abs{\alpha}\leq q_\ell,
\end{equation*}
where $m^{\beta}_{\omega'}\!(f)$ is a moment of the function $f$, 
defined by
\begin{equation*}
m^\beta_{\omega'}(f) = \int_{\omega'} 
   (\yb-\xb_{\omega'})^\beta f(\yb) \, d\yb,\quad \abs{\beta}\leq q_\ell.
\end{equation*}
Since the relationship between the moments and expansion coefficients is
linear, it is written in matrix form as
$\lambda_\omega = T(\omega,\omega') m_{\omega'}$.  The moments in a
leaf node are computed by numerical quadrature, otherwise the moments
$\omega$ can be computed by from the moments of the children. The
latter is also a linear translation written as
$m_\omega = M(\omega,\omega') m_{\omega'},\; \omega'\in
\Kids(\omega)$.  Once all expansion coefficients are computed in all
levels, they are agglomerated by translating coefficients from
the parent to the children until a leaf node is reached. The
corresponding matrix is denoted by
$L(\omega',\omega),\; \omega'\in \Kids(\omega)$. Finally, the
agglomerated series expansion in a leaf node is integrated
against the basis functions,
\begin{equation*}
  u_{\omega}^\alpha = \sum_{\abs{\beta}\leq q_\ell} \int_{\omega} \phi^\alpha_\omega(\xb)
  (\xb-\xb_\omega)^\beta \,d\xb \,\lambda_{\beta,\omega}
  , \quad \abs{\alpha}\leq p_\ell,
\end{equation*}
which in matrix notation is $u_\omega = U(\omega) \lambda_{\omega}$.
The complete procedure is summarized in algorithm \ref{algo:vtv}.

The evaluation of a nearfield interaction
$\inner{\phi^\alpha_\omega}{\N_{\omega'} f}{L_2(\Omega)}$ involve
singular integrals over Cartesian products of two tetrahedra.
For the case of singular integrals over triangles, 
there are well known transformations that
convert the singular integral to an integral over a four dimensional
cube with a smooth integrand, which in turn can be treated with tensor
product Gauss-Legendre quadrature, see~\cite{sauter-schwab11}.  For
the integrals required here similar singularity removing
transformations can be constructed, that result in smooth integrals
over six dimensional hypercubes, more details can be found in 
the PhD dissertation
\cite{mohyaddin21}.

\begin{algorithm}
\caption{Boundary concentrated FMM for the VtV calculation.}
\label{algo:vtv}
\begin{algorithmic}
\For{ $\ell=0:L$}\Comment{Nearfield Calculation.}
\For{ $\omega \in L^*_\ell$}
\State $u_{\omega,\alpha} = \sum\limits_{\omega' \in \Nbrs^*(\omega)}
\inner{\phi_\omega^\alpha}{\tilde \N_{\omega'} f}{L_2(\omega)},\;
  \abs{\alpha}\leq p_\ell$
\EndFor
\EndFor
\\
\For{ $\ell=0:L$}\Comment{Moment Calculation.} 
\For{ $\omega \in L_\ell^*$} 
\State $m_\omega^\alpha =
       \inner{(\cdot-\xb_\omega)^\alpha}{f}{L_2(\omega)},
       \;\abs{\alpha}\leq q_\ell $
\EndFor
\EndFor
\\
\For{ $\ell=L-1:0$}\Comment{Upward Pass.} 
\For{ $\omega \in M_\ell$}
\State $m_\omega = \sum\limits_{\omega' \in \Kids(\omega)}  M(\omega, \omega') m_{\omega'}$
\EndFor
\EndFor
\\
\For{ $\ell=L:0$}\Comment{Interaction Phase.}
\For{ $\omega \in C_\ell$}
\State $\lambda_\omega = \sum\limits_{\omega' \in \Inter(\omega)}  T(\omega, \omega') m_{\omega'}$
\EndFor
\EndFor
\\
\For{ $\ell=0:L-1$}\Comment{Downward Pass.} 
\For{ $\omega \in M_\ell$ and $\omega' \in \Kids(\omega)$}
\State $\lambda_{\omega'}\;+\!\!\!= L(\omega', \omega) \lambda_\omega$
\EndFor
\EndFor
\\
\For{ $\ell=0:L$}\Comment{Evaluation Phase.}
\For{$\omega \in L^*_\ell$}
\State $u_\omega\; +\!\!\!=  U(\omega) \lambda_\omega$
\EndFor
\EndFor
\end{algorithmic}
\end{algorithm}

The algorithms for the VtB and the BtV calculation are based on the
splittings \eqref{spaceDecomp:OmGa} and \eqref{spaceDecomp:GaOm}.  The
required changes for replacing the source or target by a surface are
obvious and not discussed in detail. The resulting algorithms are given
in~\ref{algo:vtb} and~\ref{algo:btv}.

Theorem~\ref{theo:cplx} implies that in all three algorithms the
number of translations is order $4^L$. The cost of each translation is
determined by the orders of basis functions $p_\ell$ and multipole
expansions $q_\ell$. In the following section we will demonstrate that
$p_\ell, q_\ell \sim L$ is sufficient to achieve convergence at any
rate implied by the regularity of the solution. Hence we have, up to
logarithmic factors, the same $\Ord(h^{-2})$ complexity as the
classical FMM for boundary to boundary calculations.

\begin{algorithm}
  \caption{Boundary concentrated FMM for the VtB calculation.}
\label{algo:vtb}
\begin{algorithmic}
\For{ $\omega \in B_L$}\Comment{Nearfield Calculation.}
\State $u_{\alpha,\omega} = \sum\limits_{\omega' \in \Nbrs(\omega)}
     \inner{\varphi^\alpha_{\gamma(\omega)}}{\N_{\omega'} f}{L_2(\gamma(\omega))},
     \;\abs{\alpha}\leq p $
\EndFor
\\
\For{ $\ell=0:L$}\Comment{Moment Calculation.}
\For{ $\omega \in L_\ell^*$} 
\State $m_\omega^\alpha =
      \inner{(\cdot-\xb_\omega)^\alpha}{f}{L_2(\omega)},
      \;\abs{\alpha}\leq q_\ell $
\EndFor
\EndFor
\\
\For{ $\ell=L-1:1$}\Comment{Upward Pass. } 
\For{ $\omega \in M_\ell$}
\State $m_\omega = \sum\limits_{\omega' \in \Kids(\omega)}  M(\omega, \omega') m_{\omega'}$
\EndFor
\EndFor
\\
\For{ $\ell=L:0$}\Comment{Interaction Phase.}
\For{ $\omega \in B_\ell$}
\State $\lambda_\omega = \sum\limits_{\omega' \in \Inter(\omega)}  T(\omega, \omega') m_{\omega'}$
\EndFor
\EndFor
\\
\For{ $\ell=0:L-1$}\Comment{Downward Pass.} 
\For{ $\omega \in B_\ell$ and $\omega' \in \Kids(\omega)\cap B_\ell$}
\State $\lambda_{\omega'}\;+\!\!\!= L(\omega', \omega) \lambda_\omega$
\EndFor
\EndFor
\\
\For{$\omega \in B_L$}\Comment{Evaluation Phase.}
\State $u_\omega\; +\!\!\!=  U(\gamma(\omega)) \lambda_\omega$
\EndFor
\end{algorithmic}
\end{algorithm}

\begin{algorithm}
  \caption{Boundary concentrated FMM for the BtV calculation.}
\label{algo:btv}
\begin{algorithmic}
\For{ $\omega \in L^*_L$}\Comment{Nearfield Calculation.}
\State $u_{\omega,\alpha} = \sum\limits_{\gamma \in \Nbrs_\Gamma(\omega)}
\inner{\phi_\omega^\alpha}{\tilde \V_{\gamma(\omega')}[\gamma_1 u] -
  \tilde \K_{\gamma(\omega')}[\gamma_0 u]}{L_2(\omega)},\;
  \abs{\alpha}\leq p_L$
\EndFor
\\
\For{ $\omega \in B_L$} \Comment{Moment Calculation.} 
\State $m_\omega^\alpha =
      \inner{(\cdot-\xb_\omega)^\alpha}{\gamma_1
        u}{L_2(\gamma(\omega))} -
      \inner{\gamma_1 (\cdot-\xb_\omega)^\alpha}{\gamma_0 u}{L_2(\gamma(\omega))},
      \;\abs{\alpha}\leq q_L $
\EndFor
\\
\For{ $\ell=L-1:0$}\Comment{Upward Pass.} 
\For{ $\omega \in B_\ell$}
\State $m_\omega = \sum\limits_{\omega' \in \Kids(\omega)\cap B_\ell}
       M(\omega, \omega') m_{\omega'}$
\EndFor
\EndFor
\\
\For{ $\ell=L:0$}\Comment{Interaction Phase.}
\For{ $\omega \in C_\ell$}
\State $\lambda_\omega = \sum\limits_{\omega' \in \Inter_\Gamma(\omega)}
   T(\omega, \omega') m_{\omega'}$
\EndFor
\EndFor
\\
\For{ $\ell=0:L-1$}\Comment{Downward Pass.}
\For{ $\omega \in M_\ell$ and $\omega' \in \Kids(\omega)$}
\State $\lambda_{\omega'}\;+\!\!\!= L(\omega', \omega) \lambda_\omega$
\EndFor
\EndFor
\\
\For{ $\ell=0:L$}\Comment{Evaluation Phase.}
\For{$\omega \in L^*_\ell$}
\State $u_\omega\; +\!\!\!=  U(\omega) \lambda_\omega$
\EndFor
\EndFor
\end{algorithmic}
\end{algorithm}

\section{Error Analysis}\label{sec:analysis:err}
The BtV and VtV algorithms are based on the assumption that layer and
volume potentials can be well approximated by a BC mesh. Since these
potentials are solutions to elliptic PDEs the error analysis in
\cite{khoromskij-melenk03} is applicable. However, since we consider
the special case of constant coefficients and an analytic source term
stronger estimates can be derived. We start with some
well known facts about Taylor series approximations of analytic
functions.

\subsection{Approximation of analytic functions by Taylor series}\label{sec:taylor}
The Taylor series of a multivariate function is obtained by
expanding the single variate function $\tau \mapsto f(\yb + \tau \hb)$
and using the chain rule. For $\tau=1$ this gives
\begin{equation*}
  f(\yb + \hb) = \sum_{n=0}^\infty D_h^n f(\yb) 
\end{equation*}
where
\begin{equation*}
  D_h^n f(\yb)  := \sum_{\abs{\alpha}=n}
  \frac{\partial^\alpha f(\yb) \hb^\alpha}{\alpha!}
  = \frac{1}{2\pi i} \int\limits_{\abs{\tau} = a}
    \frac{f(\yb + \tau \hb)}{\tau^{n+1}} \,d\tau.
\end{equation*}
The integral representation of $D_h^n f$ is a simple consequence of
Cauchy's integral formula. Assuming that $f$ is analytic in a complex
neighborhood of $\Omega$ then for $\yb\in \Omega$ and
$\abs{\hb} \leq \dist(\yb,\Gamma)$ we can set
$a=1/\abs{\hb}$. Estimating the integral in the obvious way gives
\begin{equation}\label{est:Dhn:f}
  \abs{D_h^n f(\yb)} \leq M \abs{\hb}^n,
\end{equation}
where $M$ is the maximum of $f$ in the neighborhood of $\Omega$. 

We assume that the Green's function depends only on the distance of
the source to the field point, that is, $G(\xb,\yb) = G(\abs{\xb-\yb})$.
Furthermore, we assume that $G(\cdot)$ is analytic except for the
origin and that there is a constant $C$ such that 
\begin{equation}\label{property:kernel}
\abs{G(\tau)} \leq {C\over \abs{\tau}},\quad 0 \not=\tau \in \C.
\end{equation}
The following estimates can be derived from the Cauchy integral formula
\begin{align}
  D^n_h G(\xb, \yb) &\leq C \frac{1}{\abs{\rb}}
  \left( \frac{\abs{\hb}}{\abs{\rb}}\right)^n,\\
  \frac{\partial}{\partial x_i} D^n_h G(\xb, \yb) &\leq C \frac{r_i}{\abs{\rb}^3}
  \left( \frac{\abs{\hb}}{\abs{\rb}}\right)^n,
\end{align}
where $\rb = \xb - \yb$ and $D^n_h$ can act on either the $\xb$ or the
$\yb$ variable. For more details, see Lemma 4.2 in \cite{tausch04}. 
The remainder of the truncated Taylor series of the Green's function is given by
\begin{equation*}
R_p(\xb, \yb) = \sum_{n=p+1}^\infty D_h^n G(\xb, \yb),
\end{equation*}
which can be estimated using a geometric series argument. For
$\abs{\hb} < \abs{\rb}$ we find
\begin{align}
\abs{R_p(\xb, \yb)} &\leq C {1 \over \abs{\rb}} 
                      \left(\frac{\abs{\hb}}{\abs{\rb}}\right)^{p+1},
                   \label{est:remain:G}\\
\abs{\frac{\partial}{\partial x_i} R_p(\xb, \yb)}
                    &\leq C \frac{1}{\abs{\rb}^2} 
                      \left(\frac{\abs{\hb}}{\abs{\rb}}\right)^{p+1}.
                      \label{est:remain:Gn}
\end{align}

\subsection{Approximation of potentials by BC meshes}
Since approximation error estimates involve derivatives we start by
estimating the derivatives or layer potentials in the point wise sense.

\subsubsection{Single Layer Potential} The single layer  potential can be written
as an $L_2(\Gamma)$-inner product
\begin{equation*}
  \tilde \V q(\xb) = \int_\Gamma G(\xb,\yb) q(\yb) \, ds_\yb
  = \inner{G(\xb,\cdot)}{q}{L_2(\Gamma)}.
\end{equation*}
For $\xb \in \Omega$ fixed, the kernel is a smooth function on
$\Gamma$ and thus differentiation and integration can be exchanged. In
addition, duality and the trace theorem implies that
\begin{equation*}
  \abs{D_h^n \tilde \V q(\xb)} \leq
  \norm{D_h^n G(\xb,\cdot)}{H^\half(\Gamma)} \norm{q}{H^{-\half}(\Gamma)}\leq
  \norm{D_h^n G(\xb,\cdot)}{H^1(\Omega^c)} \norm{q}{H^{-\half}(\Gamma)}.
\end{equation*}
In the last step we applied the trace theorem to the exterior domain
$\Omega^c = \R^3\setminus \bar \Omega$ to avoid the singularity of the
Green's function. To estimate the last term let
$r = \dist(\xb,\Gamma)$ and $B_r^c(\xb)$ the exterior of the sphere of
radius $r$ centered in $\xb$. Since $\Omega^c \subset B_r^c(\xb)$ we
have with the estimates of section~\ref{sec:taylor}
\begin{align*}
\norm{D_h^n G(\xb,\cdot)}{H^1(\Omega^c)}^2
  &\leq \norm{\nabla D_h^n G(\xb,\cdot)}{L_2(B^c_r(\xb))}^2  +
    \norm{D_h^n G(\xb,\cdot)}{L_2(B^c_r(\xb))}^2 \\
  &\leq C \int_{r}^\infty \left( \frac{1}{\rho^{2n+4}} +
    \frac{1}{\rho^{2n+2}} \right)\rho^2 \,d\rho \abs{\hb}^{2n} \\
  &\leq C  \frac{1}{r} \left( \frac{\abs{\hb}}{r}\right)^{2n}.
\end{align*}
Note that we have absorbed a negative power of $n$ in the
constant, as it is not relevant in the following error estimate. Thus
\begin{equation}\label{est:Dhn:V}
  \abs{D_h^n \tilde \V q(\xb)} \leq C  \frac{1}{r^\half}
  \left( \frac{\abs{\hb}}{r}\right)^{n} \norm{q}{H^{-\half}(\Gamma)}.
\end{equation}

\subsubsection{Double Layer Potential} This estimate follows along similar lines,
\begin{align*}
  \abs{D_h^n \tilde \K u(\xb)} &= \abs{\inner{\gamma_1 D_h^n G(\xb,\cdot)}{u}{L_2(\Gamma)}}\\
  &\leq \norm{D_h^n \gamma_1 G(\xb,\cdot)}{H^{-\half}(\Gamma)} \norm{u}{H^{\half}(\Gamma)}
    \leq \norm{D_h^n G(\xb,\cdot)}{H^1(\Omega^c)} \norm{u}{H^{\half}(\Gamma)}.
\end{align*}
Here we have used that $G(\xb,\cdot)$ is in the kernel of
$\mathcal{L}$ and that $\gamma_1 : H^1(\Omega^c, \mathcal{L}) \to
H^{-\half}(\Gamma)$ is continuous. This implies that 
\begin{equation}\label{est:Dhn:K}
  \abs{D_h^n \tilde \K u(\xb)} \leq C  \frac{1}{r^\half}
  \left( \frac{\abs{\hb}}{r}\right)^{n} \norm{u}{H^{\half}(\Gamma)}.
\end{equation}

\subsubsection{Estimates for Higher Regularity} 
Estimates \eqref{est:Dhn:V} and  \eqref{est:Dhn:K}
can be improved if the potential has more than $H^1(\Omega)$
regularity. To that end, suppose $u$ is either the single or double layer potential
with a source term that generates a potential with $H^{s+1}(\Omega)$ regularity for some integer
$s \geq 1$. 
Then $D_h^s u \in H^1(\Omega)$ 
is a solution of $\mathcal{L} D_h^s u = 0$. Applying the Green's
representation formula \eqref{green:ie:vol} to $D_h^s u$ implies that
\begin{equation*}
  D_h^s u = \tilde \V [\gamma_1 D_h^s u] - \tilde \K [\gamma_0 D_h^s u],
  \qquad \text{in}\;\Omega.
\end{equation*}
Now use estimates \eqref{est:Dhn:V}  and \eqref{est:Dhn:K}
for the potentials on the right hand side, and the fact that
$\norm{\gamma_1 D_h^s u}{H^{-\half}(\Gamma)}$ and 
$\norm{\gamma_0 D_h^s u}{H^{\half}(\Gamma)}$ are bounded by 
$c \norm{D_h^s u}{H^1(\Omega)}$. Thus 
\begin{equation}\label{est:Dhn:us}
\begin{aligned}
  \abs{D_h^n u(\xb)} &= \frac{(n-s)! s!}{n!} \abs{D_h^{n-s} D_h^s u(\xb)}
  \leq 
  C  \frac{1}{r^\half}
  \left( \frac{\abs{\hb}}{r}\right)^{n-s}\norm{D_h^s u}{H^{1}(\Omega)}\\
  &\leq C  \frac{\abs{\hb}^s}{r^\half} \left( \frac{\abs{\hb}}{r}\right)^{n-s}
 \norm{u}{H^{s+1}(\Omega)}.
\end{aligned}
\end{equation}

\subsubsection{Volume Potential} The difficulty of estimating
derivatives of the volume potential is that the evaluation point is
inside the domain of integration and higher derivatives of the kernel
are strongly singular. This issue can be resolved with integration by
parts. Since $\partial_{i,\xb} G(\xb,\cdot)$ is still
weakly singular we have
\begin{align*}
\partial_i \tilde \N f(\xb) &= \int_\Omega \partial_{i,\xb} G(\xb,\yb) f(\yb) d\yb
= - \int_\Omega \partial_{i,\yb} G(\xb,\yb) f(\yb) d\yb\\
&= \int_\Omega G(\xb,\yb) \partial_i f(\yb) d\yb
  - \int_\Gamma G(\xb,\yb) f(\yb) n_{i,\yb} ds_\yb
\end{align*}
Now the right hand side can be differentiated to obtain second
derivatives of $\tilde \N f(\xb)$. 
Repeated differentiation and partial integration results in the following
formula
\begin{align*}
  D_h^n \tilde \N f(\xb) &=
  \int_\Omega G(\xb,\yb) D_h^n f(\yb) \,d\yb \\
  &+ \sum_{s=1}^{n} \frac{(n-s)! (s-1)!}{n!} 
    \int_\Gamma   D_{\xb,h}^{n-s} G(\xb,\yb) D_h^{s-1} f(\yb)
    \, \hb \cdot\mathbf{n}_\yb   \, ds_\yb.
\end{align*}
With the estimates of section~\ref{sec:taylor} it follows that the
volume integral has the upper bound $MC \abs{\hb}^n$ and is therefore
of lower order. The boundary terms can be estimated in a similar
fashion as the single layer potential. Here we note that
$\norm{D_h^{s-1} f}{H^{-\half}(\Gamma)} \leq \norm{D_h^{s-1}
  f}{L_2(\Gamma)} \leq M \abs{\hb}^{s-1}$. This leads to
\begin{equation*}
  \abs{D_h^n \tilde \N f(\xb)} \leq CM \abs{\hb}^n + 
  \frac{CM}{r^\half}
  \left(   \frac{\abs{\hb}}{r} \right)^n
  \sum_{s=1}^{n}  \frac{(n-s)! (s-1)!}{n!} r^s . 
\end{equation*}
Since the fractions in the sum are bounded by unity, and we only
consider small $r$ the sum can be
bounded by a factor $Cr$. It follows that
\begin{equation}\label{est:Dhn:N}
  \abs{D_h^n \tilde \N f(\xb)} \leq CM r^\half
  \left( \frac{\abs{\hb}}{r} \right)^n.
\end{equation}
If $f$ is such that the volume potential $u$ has
$H^{s+1}(\Omega)$-regularity, then we can apply the same trick as
before. Since $\mathcal{L} D_h^s u =   D_h^s f$ and  $D_h^s u \in
H^1(\Omega)$ we have from Green's representation formula that
\begin{equation*}
  D_h^s u = \tilde \V [\gamma_1 D_h^s u] - \tilde \K [\gamma_0 D_h^s u]
         + \tilde \N D_h^s f
  \qquad \text{in}\;\Omega.
\end{equation*}
Applying estimates \eqref{est:Dhn:us} and \eqref{est:Dhn:N} to the
potentials in the last equation leads to the estimate
\begin{equation}
  \abs{D_h^n u(\xb)} \leq C M_s \frac{\abs{\hb}^s}{r^\half}
  \left( \frac{\abs{\hb}}{r} \right)^{n-s}
\end{equation}
where
\begin{equation*}
M_s = \max \left\{ M, \norm{u}{H^{s+1}} \right\}.
\end{equation*}

\subsubsection{Approximation Errors}
We now consider the error of best approximation in the $L_2(\Omega)$-norm
if a potential is approximated by a function in the finite
element space $S^\Omega_h$ defined in \eqref{def:femspace:volume}.
We first consider all tetrahedra in the sets $L_\ell$,
which are well separated from the boundary.
Here we can approximate the potential $u$ (which can either be a
surface or volume potential) by a
$p$-th order Taylor series centered in the center $\xb_\omega$ of the tetrahedron. If 
$\hb = \xb - \xb_\omega$ and $r = \dist(\xb_\omega,\Gamma)$ then
\begin{equation*}
  \abs{u(\xb) - T_p u(\xb)} \leq
  \sum_{n=p+1}^\infty \abs{D_h^n u(\xb_\omega)} 
  \leq C  \frac{\abs{\hb}^s}{r^\half} \sum_{n=p+1}^\infty \left( \frac{\abs{\hb}}{r}\right)^{n-s}
 \norm{u}{H^{s+1}(\Omega)}
\end{equation*}
If $\xb_\Gamma \in \omega' \in B_\ell$ is the closest point of $\xb_\omega$ on
$\Gamma$, then it follows from the definitions \eqref{def:diam} and
\eqref{def:sepRatio} that
\begin{equation*}
  \frac{\abs{\hb}}{r} \leq
  \frac{\rho_\omega}{\abs{\xb_\omega - \xb_{\omega'} + \xb_{\omega'} - \xb_\Gamma}} \leq
  \frac{\rho_\omega}{\abs{\xb_\omega - \xb_{\omega'}} - \rho_{\omega'}} \leq
  \frac{\rho_\omega + \rho_{\omega'}}{\abs{\xb_\omega - \xb_{\omega'}}} =
  \eta(\omega,\omega') \leq \eta
\end{equation*}
where in the last step we used that $\omega' \not\in \Nbrs(\omega)$. 
Since $r \geq c 2^{-\ell}$, $\abs{\hb} \leq c 2^{-\ell}$ and $\eta<1$, we obtain from the
geometric series that
\begin{equation*}
  \abs{u(\xb) - T_p u(\xb)} \leq C \eta^p  2^{-\ell(s-\half)} \norm{u}{H^{s+1}(\Omega)}.
\end{equation*}
Since $s$ is fixed and small the constant absorbs a factor of $\eta^{-s}$.
For the $L_2$-orthogonal projector $P_\omega$ into the subspace of degree-$p$
polynomials we get 
\begin{equation*}
  \norm{u - P_\omega u}{L_2(\omega)}^2 \leq \norm{u - T_p
    u}{L_2(\omega)}^2 \leq
  \abs{\omega} \max_{\xb\in\omega} \abs{u(\xb) - T u(\xb)}^2 
  \leq C \eta^{2p}  2^{-2\ell(s+2)},
\end{equation*}
since $\abs{\omega} \leq c 2^{-3\ell}$.  For the remaining
$\omega \in M_L$, the standard estimate for $P_\omega$ can be
applied. Thus
\begin{equation*}
  \norm{u - P_\omega u}{L_2(\omega)} \leq C \abs{\hb}^{p+1} \norm{u}{H^{p+1}(\omega)},
  \quad 0\leq p \leq s.
\end{equation*}
Summing over all $\omega \in \Triang$ gives the total error. Since the number of
tetrahedra in $L_\ell$ is bounded by $c 2^{2\ell}$ we get 
\begin{align*}
  \norm{u - P_\Triang u}{L_2(\Omega)}^2 &\leq
  \sum_{\ell=0}^L \sum_{\omega \in L_\ell} \norm{u - P_\omega u}{L_2(\omega)}^2
   + \sum_{\omega \in M_L} \norm{u - P_\omega u}{L_2(\omega)}^2\\
  &\leq C \left(
  \sum_{\ell=0}^L \eta^{2p_\ell} 2^{-2 \ell s} 
    + \abs{\hb}^{2p+2}\right) \norm{u}{H^{s+1}(\Omega)}^2
\end{align*}
The last estimate makes clear how the choice of $\eta$ and the expansion order $p_\ell$
for tetrahedra in $L_\ell$ and the order $p$ for tetrahedra in $M_L$ affect
the accuracy. In particular, if we let $p_\ell = L-\ell$, then
\begin{equation}\label{err:approx}
  \sum_{\ell=0}^L \eta^{2p_\ell} 2^{-2 \ell s} =
  2^{-2 L s} \sum_{\ell=0}^L (\eta 2^s)^{2\ell}.
\end{equation}
In order to bound the geometric series, $\eta$ must satisfy
$\eta < 2^{-s}$. Moreover, the order $p$ for the tetrahedra
in $M_L$ must satisfy $p=s-1$. This leads to the error
\begin{equation}\label{est:approx:err}
  \norm{u - P_\Triang u}{L_2(\Omega)} \leq C h^s  \norm{u}{H^{s+1}(\Omega)}.
\end{equation}

\subsection{FMM errors}
To analyze the VtV algorithm we return to the space decomposition
\eqref{spaceDecomp:OmOm}, which implies the following splitting of the
bilinear form induced by the volume potential
\begin{align*}
  \inner{v_h}{\tilde \N f}{L_2(\Omega)} &=
  \sum_{\ell=0}^L \sum_{\omega \in L^*_\ell\atop \omega'\in \Nbrs^*(\omega)} 
    \int_{\omega}\int_{\omega'} G(\xb,\yb) v_h(\xb) f(\yb) \,d\yb d\xb\\
    &+  \sum_{\ell=0}^L \sum_{\omega \in C_\ell\atop \omega'\in \Inter(\omega)}
    \int_{\omega}\int_{\omega'} G(\xb,\yb) v_h(\xb) f(\yb) \,d\yb d\xb.
\end{align*}
In the fast algorithm the integrals in the first sum are computed
directly, whereas in the second sum the Green's function is replaced
by its Taylor series approximation. Denoting the approximate volume potential by
$\tilde \N_h$, we get
\begin{equation*}
  \inner{v_h}{(\tilde \N - \tilde \N_h) f}{L_2(\Omega)} =
\sum_{\ell=0}^L \sum_{\omega \in C_\ell\atop \omega'\in \Inter(\omega)}
    \int_{\omega}\int_{\omega'} R_{q_\ell}(\xb,\yb) v_h(\xb) f(\yb) \,d\yb d\xb.
\end{equation*}
where $q_\ell$ is the expansion order, which may depend on the level,
and $R_{q_\ell}$ is the remainder of the Taylor series. In the
estimates for the remainder \eqref{est:remain:G}, we set
$\rb=\rb_{\omega,\omega'}=\xb_\omega - \xb_{\omega'}$ and
$\hb = \xb-\xb_\omega - \yb + \xb_{\omega'}$. Then it follows that
$\abs{\hb}/\abs{\rb} \leq \eta(\omega,\omega')\leq \eta$, because
$\omega$ and $\omega'$ are in interaction lists. Thus we can estimate
for $v_h \in S_h^\Omega$,
\begin{align*}
  \abs{\inner{v_h}{(\tilde \N - \tilde \N_h) f}{L_2(\Omega)}} &\leq
C \sum_{\ell=0}^L \sum_{\omega \in C_\ell\atop \omega'\in \Inter(\omega)}
\frac{1}{\abs{\rb_{\omega,\omega'}}} \eta^{q_\ell}
\int_{\omega}\abs{v_h(\xb)}\, d\xb \int_{\omega'} \abs{f(\yb)} \,d\yb,\\
&\leq C
\sum_{\ell=0}^L \sum_{\omega \in C_\ell\atop \omega'\in \Inter(\omega)}
\frac{\abs{\omega}^\half\abs{\omega'}^\half }{\abs{\rb_{\omega,\omega'}}}     
\eta^{q_\ell} \norm{v_h}{L_2(\omega)}\norm{f}{L_2(\omega')}.
\end{align*}
To continue the estimate note that the number of terms in an
interaction list is uniformly bounded, hence 
\begin{equation*}
\sum_{\omega \in C_\ell\atop \omega'\in \Inter(\omega)}
  \norm{v_h}{L_2(\omega)} \norm{f}{L_2(\omega')} \leq
C \norm{v_h}{L_2(\Omega)} \norm{f}{L_2(\Omega)}.
\end{equation*}
Moreover, because of the shape regularity we have
$\abs{\omega}, \abs{\omega}' \leq c 2^{-3\ell}$ and
$\abs{\rb_{\omega,\omega'}} \geq c 2^{-\ell}$. Adding the contribution
of each level $\ell$ gives the estimate
\begin{equation}\label{est:bilinear:eps}
  \abs{\inner{v_h}{(\tilde \N - \tilde \N_h) f}{L_2(\Omega)}}  \leq
C \epsilon_2
\norm{v_h}{L_2(\Omega)}\norm{f}{L_2(\Omega)}
\end{equation}
where
\begin{equation*}
\epsilon_r = \sum_{\ell=0}^L 2^{-r\ell} \eta^{q_\ell} .
\end{equation*}
For the VtV calculation we have $r=2$, but considering a general
value of $r$ will facilitate the discussion for the VtB and BtV algorithms.
The goal is to determine the expansion orders $q_\ell$ such that $\epsilon_r$
is of the same magnitude as the approximation error in
\eqref{err:approx}. Note that we use the same $\eta<2^{-s}$. Motivated
by the discussion of the approximation error, we consider expansion
orders dependent on the level as follows
\begin{equation*}
  q_\ell = q_0 + L - \ell
\end{equation*}  
then
\begin{equation*}
  \epsilon_r = \eta^{q_0}  \sum_{\ell=0}^L 2^{-r\ell} \eta^{L-\ell}
  \leq \eta^{q_0} \max\left\{ 2^{-r}, \eta \right\}^L
  \leq 2^{-s q_0}  \max\left\{ 2^{-r}, 2^{-s}\right\}^L.
\end{equation*}  
Hence, if the expansion order in the finest level is given by
\begin{equation}\label{def:q0}
  q_0 = \begin{cases}
    \frac{s-r}{s} L &\text{if}\; r<s\\
    0 &\text{if}\; r\geq s
  \end{cases}
\end{equation}  
then it follows that $\epsilon_r \leq C 2^{-sL}$ and thus
\begin{equation}\label{est:bilinear}
  \abs{\inner{v_h}{(\tilde \N - \tilde \N_h) f}{L_2(\Omega)}}  \leq
C 2^{-sL} \norm{v_h}{L_2(\Omega)}\norm{f}{L_2(\Omega)}.
\end{equation}  
The overall error can be obtained with a Strang-Lemma type
argument. To that end, note that
$\inner{v}{u_h}{L_2(\Omega)} = \inner{P_\Triang v}{u_h}{L_2(\Omega}$
holds for all $u_h \in S_h^\Omega$ and $v \in L_2(\Omega)$. Moreover,
$\norm{P_\Triang v}{L_2(\Omega)} \leq \norm{v}{L_2(\Omega)}$. Then for
$u_h = P_\Triang \tilde \N f$ and $u_h^f = \tilde \N_h f$ we obtain
\begin{align*}
  \norm{u_h - u_h^f}{L_2(\Omega)} &= \sup_{v\in L_2(\Omega)} 
  \frac{\inner{v}{u_h - u_h^f}{L_2(\Omega)}}{\norm{v}{L_2(\Omega)}}
  \leq \sup_{v\in L_2(\Omega)} 
  \frac{\inner{P_\Triang v}{u_h - u_h^f}{L_2(\Omega)}}
        {\norm{P_\Triang v}{L_2(\Omega)}}  \\
&\leq \sup_{v_h\in S_h} 
\frac{\inner{v_h}{u_h - u_h^f}{L_2(\Omega)}}{\norm{v_h}{L_2(\Omega)}}
\leq 2^{-sL} \norm{f}{L_2(\Omega)}.
\end{align*}
where in the last step we used estimate \eqref{est:bilinear}. The
overall error involves the approximation error \eqref{est:approx:err}
and the triangle inequality. We get
\begin{equation*}
  \norm{u - u_h^f}{L_2(\Omega)} \leq
  \norm{u_h - P_\Triang u_}{L_2(\Omega)} + \norm{P_\Triang u - u_h^f}{L_2(\Omega)} 
  \leq C h^s  \norm{u}{H^{s+1}(\Omega)}.
\end{equation*}

The error analysis for the remaining potential calculations is
completely analogous. Revisiting the calculations that led to
\eqref{est:bilinear:eps} makes clear that
\begin{align*}
  \abs{\inner{v_h}{(\tilde \V - \tilde \V_h) q}{L_2(\Omega)}}  &\leq
C \epsilon_{\thrf}
\norm{v_h}{L_2(\Omega)}\norm{q}{L_2(\Gamma)}, 
  \quad v_h \in S^{\Omega}_h,\; q\in L_2(\Gamma), \\
  \abs{\inner{v_h}{(\tilde \K - \tilde \K_h) q}{L_2(\Omega)}}  &\leq
C \epsilon_{\half}
\norm{v_h}{L_2(\Omega)}\norm{q}{L_2(\Gamma)}, 
  \quad v_h \in S^{\Omega}_h,\; q\in L_2(\Gamma), \\
\abs{\inner{w_h}{(\N - \N_h) f}{L_2(\Gamma)}}  &\leq
C \epsilon_{\thrf}
\norm{w_h}{L_2(\Gamma)}\norm{f}{L_2(\Omega)}, 
  \quad w_h \in S^{\Gamma}_h,\; f\in L_2(\Omega).
\end{align*}
Note that the reduced values of $r$ come from the fact that one of the
functions is defined on a surface. Moreover, when the double
layer potential is evaluated in the volume, the kernel produces a
stronger singularity than the single layer potential which results in
a further reduction of the value of $r$. However, if the expansion
order in the finest level is determined as in \eqref{def:q0} it is
still possible to obtain $2^{-sL}$ convergence of the overall error.

\section{Numerical Results}\label{sec:numresults}
We first give some details about the meshes used in our numerical experiments.
The domain is the cube $\Omega = [-1,1]^3$, where the level-zero
refinement consists of 48 congruent tetrahedra. In all experiments we
use $\eta_0=0.5$ in the definition of neighbors in
\eqref{def:neighbors}. The data for the resulting BC mesh is displayed
in table~\ref{tab:bcRefinement}. Here, $\Nbrs_{max}$ and $\Inter_{max}$
are the maximal number of neighbors and interaction lists in the level
and CPU is the CPU time in seconds to
compute the mesh and neighbor and interaction lists up to the given
refinement. 

It is apparent that the first three
refinements are uniform and thus the first leaves show up in level
three. Beginning with level six one can see that the asymptotic
estimates of theorem~\ref{theo:cplx} are reproduced. The maximal number of
neighbors in a level converges to a perhaps unexpectedly large
number, but this can be explained by the fact that we fill a three
dimensional domain with tetrahedra. Our implementation stores the
moments and expansion coefficients as well as the neighbor and
interaction lists for each tetrahedron. The latter turns out to be the
dominant memory usage.

\begin{table}
\begin{center}
\begin{tabular}{c c c c c c c c c}
  $\ell$  &  $\# C_{\ell}$ & fac & $\# L_{\ell}$& fac
  & $\Nbrs_{max}$ & $\Inter_{max}$ & CPU & fac\\
    \hline
    $0$  &  48     &       & 0        & & 48   & 0   & &\\
    $1$  & 384     & 8     & 0        & & 380  &178  & &\\
    $2$  & 3072    & 8     & 0        & & 967  &2249 & &\\
    $3$  & 24576   & 8     & 3552     & & 1273 &6872 & & \\
    $4$  & 168192  & 6.8   & 58272    & 16.4 & 1302 &8120 & 111\\
    $5$  & 879360  & 5.2   & 379680   & 6.5 & 1302 &8120 & 643  & 5.8\\
    $6$  & 3997440 & 4.5   & 1870368  & 4.9 & 1302 &8120 & 3121 & 4.9\\
\end{tabular}
\caption{Mesh data.}
\end{center}
\label{tab:bcRefinement}  
\end{table}
In all numerical experiments reported below, we use $p=0$ and $p_\ell
= L-\ell$ for the expansion order of the tetrahedra. Thus according to
\eqref{est:approx:err} we expect that the best approximation error of the BC
mesh is $O(h)$. For the multipole expansion orders in \eqref{def:q0}
we have set $q_\ell = 4 + L-l$. 

\begin{figure}
\begin{center}
\includegraphics[width=0.4\textwidth]{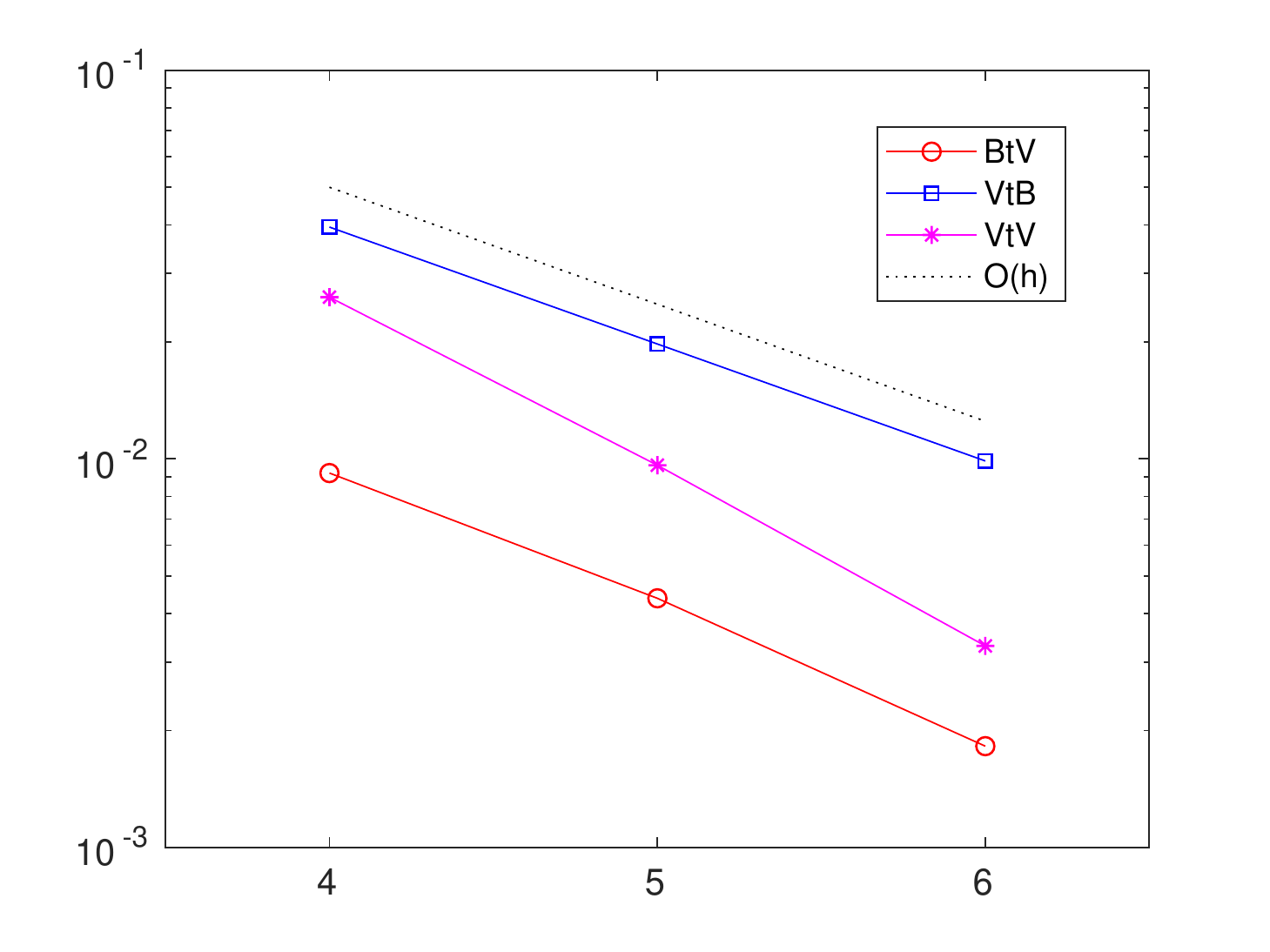}
\includegraphics[width=0.4\textwidth]{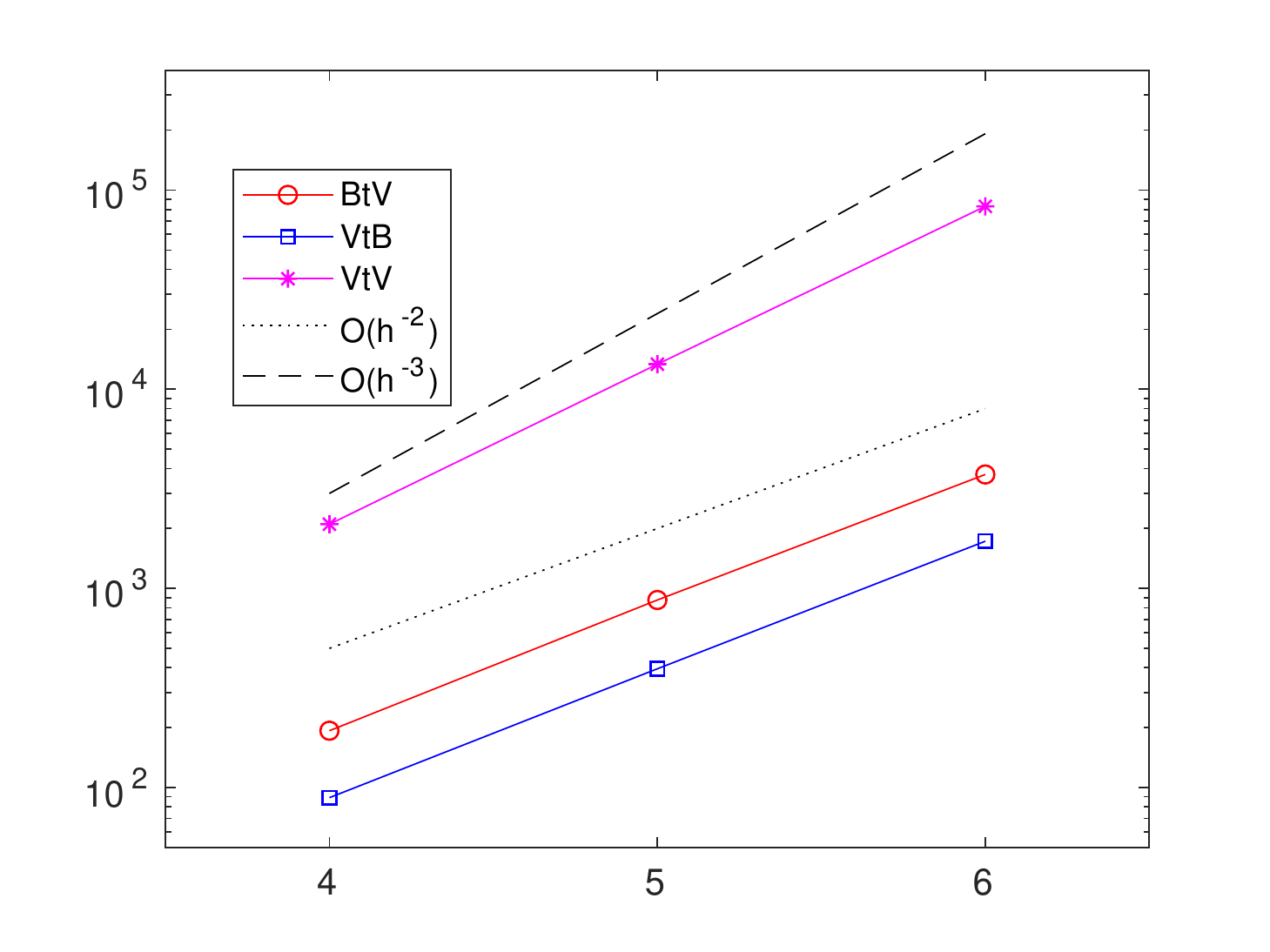}
\end{center}
\caption{
  Errors and CPU times for the potential calculation versus the mesh
  refinement level.} 
\label{fig:Results}
\end{figure}

We now illustrate the behavior of the above algorithms on examples with
known potentials. The errors and CPU timings are displayed in \ref{fig:Results}.
In particular, we consider the functions
\begin{align*}
  u_L(x,y,z) &= \left( (x-x_0)^2 + (y-y_0)^2 + (z-z_0)^2 \right)^{-\half},\\
  u_P(x,y,z) &= \exp\left( -r^2 \right),\\
    f(x,y,z) &= \left(6 - 4 r^2\right)\exp\left(-r^2\right),   
\end{align*}
where $(x_0, y_0, z_0)$ = $(3,0,0)$ and
$r = (x^2 + y^2 + z^2)^\half$. The function $u_L$ solves the Laplace
equation and the function $u_P$ solves the Poisson equation $-\Delta u_P
= f$.

To test the BtV algorithm we compute the right hand side in the
Green's representation formula
\begin{equation*}
  u_L = \tilde \V[\gamma_1 u_L] - \tilde \K [\gamma_0 u_L] 
\end{equation*}  
and compare the calculated potential with the analytic value on the
left hand side of the equation. 
To test the VtB algorithm we compute the right hand side in the
Green's integral equation
\begin{equation*}
  \gamma_0 u_P = 2\left\{ \V[\gamma_1 u_P] - \K [\gamma_0 u_P] + \N f \right\}
\end{equation*}  
and compare the calculated potential with the analytic value on the
left hand side of the equation. The evaluation of the right hand side
also involves a BtB calculation, which can be performed with the
algorithm that is obtained by restricting the source and target
domains of the VtV calculation on the boundary. The result 
is the standard FMM for surface potentials, see,
e.g.,~\cite{tausch04}. In \ref{fig:Results} we report the combined
error and the CPU times for evaluating $\N f$.

Finally, to test the VtV algorithm, we compute the right hand side in the
Green's representation formula for the Poisson equation
\begin{equation*}
  u_P = \tilde \V[\gamma_1 u_P] - \tilde \K [\gamma_0 u_P] + \tilde \N f
\end{equation*}  
and compare the calculated potential with the analytic value on the
left hand side of the equation. The evaluation of the right hand side
also involves a BtV calculation, which we have already tested.
Even though $u_P$ is in $C^\infty(\bar \Omega)$, it
can be expected that the individual potentials have much lower
regularity in the domain. Since the different potential calculations
work independently this indicates that the individual algorithms also work
for lower regularity situations.

The order of
the Gauss-Legendre rule for the nearfield critically influences the cost and
accuracy of the overall algorithm. In our implementation, we use a
fixed quadrature order in the finest level. For the coarser level
nearfield interactions of the VtV algorithm,
the order is increased in each level.

As it is apparent from figure \ref{fig:Results} the errors of all
potential calculations converge at
the expected $O(h)$ rate, the VtV result appears 
faster, probably because the multipole error in
\eqref{est:bilinear:eps} gives smaller estimates
for volumes. The timing of the VtB and BtV
methods are in excellent agreement with the theoretical $O(h^{-2})$
estimate. The data for the VtV algorithm is somewhat higher, but
considerably better than $O(h^{-3})$. A likely
cause is that this algorithm evaluates nearfield interactions 
in coarser levels. However, in table~\ref{tab:bcRefinement}
it is apparent that the number of leaves behaves
pre-asymptotically in the coarser levels, and therefore the calculated
number of levels are not yet sufficient exhibit the expected $O(h^{-2})$ complexity.


\end{document}